\documentclass[11pt]{article}

\usepackage{amsmath,amsthm}

\usepackage{amssymb,latexsym}

\usepackage{enumerate}

\newtheorem{tw}{Theorem}
\newtheorem{corollary}[tw]{Corollary}
\newenvironment{dow}{\it Proof.\rm}{\hfill $\Box$}

\theoremstyle{definition}
\newtheorem*{definition}{Definition}
\newtheorem*{uw}{Remark}

\newcommand{\BR}{{\mathbb R}}
\newcommand{\BX}{{\mathbb X}}

\newcommand{\BB}{{\mathcal{B}}}
\newcommand{\FF}{{\mathcal{F}}}
\newcommand{\LL}{{\mathcal{L}}}

\begin{document}

\title {The early exercise premium representation for
American options on multiply assets}
\author {Tomasz Klimsiak and Andrzej Rozkosz
\mbox{}\\[2mm]
{\small  Faculty of Mathematics and Computer Science,
Nicolaus Copernicus University},\\
{\small Chopina 12/18, 87-100 Toru\'n, Poland}}
\date{}
\maketitle

\begin{abstract}
\noindent {\bf Abstract} In the paper we  consider the problem of
valuation of American options written on dividend-paying assets
whose price dynamics follow the classical multidimensional Black
and Scholes model. We provide a general early exercise premium
representation formula for options with payoff functions which are
convex or satisfy mild regularity assumptions. Examples include
index options, spread options, call on max options, put on min
options, multiply strike options and power-product options. In the
proof of the formula we exploit close connections between the
optimal stopping problems associated with valuation of American
options, obstacle problems and reflected backward stochastic
differential equations.
\end{abstract}


\footnotetext{Email addresses: tomas@mat.umk.pl (T. Klimsiak),
rozkosz@mat.umk.pl (A. Rozkosz). }

\footnotetext{Research supported by National Science Centre grant
no. 2012/07/B/ST1/03508.}

\section{Introduction}

In the paper we study American options written on dividend-paying
assets. We assume that the underlying assets dynamics follow the
classical multidimensional Black and Scholes model. It is now well
known that the arbitrage-free value  of American options can be
expressed in terms of the optimal stopping problem (Bensoussan
\cite{Be}, Karatzas \cite{K}; see also \cite{KS} for nice
exposition and additional references), in terms of variational
inequalities (Jaillet, Lamberton and Lapeyre \cite{JLL}) and in
terms of solutions of reflected BSDEs (El Karoui and Quenez
\cite{EQ}). Although these approaches provide complete
characterization of the option value (see Section \ref{sec2} for a
short review), the paper by Broadie and Detemple \cite{BD2} shows
that it is of interest to provide alternative representation,
which expresses the value of an American option as the value of
the corresponding European option plus the gain from early
exercise. The main reason is that the representation, called the
early exercise premium formula, gives useful information on the
determinants of the option value. The formula was proved first by
Kim \cite{Ki} in the case of standard American put option on a
single asset. Another important contributions in the case of
single asset include Broadie and Detemple \cite{BD1}, El Karoui
and Karatzas \cite{EK} and Jacka \cite{J} (see also \cite{D,KS}
and the references therein). The case of options on multiply
dividend-paying  assets is more difficult and has received rather
little attention in the literature. In the important paper
\cite{BD2} and next in Detemple, Feng and Tian \cite{DFT} (see
also \cite{D}) the early exercise premium formula was established
for concrete classes of options on multiply assets. Note that in
the last paper call on min option, i.e. option with nonconvex
payoff function is investigated. A subclass of call on min options
consisting of capped options is studied in \cite{BD1,BD2,BD3} (see
also \cite{D}).

In the present paper we provide a unified way of treating a wide
variety of seemingly disparate examples. It allows us to prove a
general exercise premium formula for options with convex payoff
functions satisfying the polynomial growth condition or payoff
function satisfying quite general condition considered in Laurence
and Salsa \cite{LS}. Verifying the last condition requires
knowledge of the payoff function and the structure of the exercise
set. Therefore it is a complicated task in general. Fortunately,
in most interesting cases one can easily check convexity of the
payoff function or check some simpler condition implying the
general condition from \cite{LS}. The class of options covered by
our formula includes index options, spread options, call on max
options, put on min options, multiply strike options,
power-product options and others.

In the proof of the exercise premium formula we rely on some
results on reflected BSDEs and their links with optimal stopping
problems (see \cite{EQ}) and with parabolic variational
inequalities established in Bally, Caballero, Fernandez and El
Karoui \cite{BCFE}. We also use classical results on regularity of
the solution of the Cauchy problem for parabolic operator with
constant coefficients, and in case of convex payoffs, some fine
properties of convex functions. Perhaps it is worth mentioning
that we do not use any regularity results on the free boundary
problem for an American option. The basic idea of the proof comes
from our earlier paper \cite{KR:BPAN} devoted to standard American
call and put options on single asset.

\section{Preliminaries}
\label{sec2}

We will assume that under the risk-neutral probability measure the
underlying assets prices $X^{s,x,1},\dots,X^{s,x,n}$ evolve on the
time interval $[s,T]$ according to stochastic differential
equation of the form
\begin{equation}
\label{eq1.1}
X^{s,x,i}_{t}=x_i+\int_{s}^{t}(r-d_{i})X^{s,x,i}_{\theta}\,d\theta
+\sum_{j=1}^n\int_{s}^{t}\sigma_{ij}X^{s,x,i}_{\theta}\,dW^j_{\theta},
\quad t\in[s,T].
\end{equation}
Here $W=(W^1,\dots,W^n)$ is a standard $n$-dimensional Wiener
process, $r\ge0$ is the rate of interest, $d_i\ge0$ is the
dividend rate of the asset $i$ and $\sigma=\{\sigma_{ij}\}$ is the
$n$-dimensional volatility matrix. We assume that
$a=\sigma\cdot\sigma^*$ is positive definite. Since the
distributions of the processes $X^{s,x,i}$ depend $\sigma$  only
through $a$, we may and will assume that $\sigma$ is a symmetric
square root of $a$. As for the payoff function $\psi$ we will
assume that it satisfies the assumptions:
\begin{enumerate}
\item[\rm{(A1)}] $\psi$ is a nonnegative continuous function on $\BR^n$
with polynomial growth,

\item[{\rm{(A2)}}]For every $t\in (0,T)$, $\psi$ is a smooth function
on $\{\psi=u\}\cap\bar Q_{t}$, i.e. there exists an open set
$U\subset \BR^n$ such that $\{u=\psi\}\cap\bar Q_{t}\subset
[0,t]\times U$ and $\psi$ is smooth on $U$ (Here
$Q_t\equiv[0,t)\times\BR^n$, $\bar Q_t\equiv[0,t]\times\BR^n$ and
$u$ is the value of an option with payoff $\psi$; see
(\ref{eq1.2}) and (\ref{eq2.4}) below)
\end{enumerate}
or
\begin{enumerate}
\item[\rm{(A3)}]$\psi$ is a nonnegative convex function on $\BR^n$
with polynomial growth.
\end{enumerate}

Note that  convex functions are locally Lipschitz, so assumption
(A3) implies (A1). Assumption (A2) is considered in \cite{LS}. It
is satisfied for instance if
\begin{enumerate}
\item[\rm{(A2$'$)}] The region where $\psi$ is strictly positive
is the union of several connected components in which $\psi$ is
smooth.
\end{enumerate}
Following \cite{LS} let us also note that unlike (A2$'$) or (A3),
condition (A2) cannot be verified by appealing to the structure of
the payoff alone. Verifying (A2) requires additional knowledge of
the structure of the exercise set  $\{u=\psi\}$.

Let $\Omega=C([0,T];\BR^n)$ and let $X$ be the canonical process
on $\Omega$. For $(s,x)\in Q_{T}$ let $P_{s,x}$  denote the law of
the process $X^{s,x}=(X^{s,x,1},\dots,X^{s,x,n})$ defined by
(\ref{eq1.1}) and let $\{\FF^s_t\}$  denote the completion of
$\sigma(X_{\theta};\theta\in[s,t])$ with respect to the family
$\{P_{s,\mu};\mu$ a finite measure on $\BB(\BR^n)\}$, where
$P_{s,\mu}(\cdot)=\int_{\BR^n}P_{s,x}(\cdot)\,\mu(dx)$. Then for
each $s\in[0,T)$, $\BX=(\Omega,(\FF^s_t)_{t\in[s,T]},X,P_{s,x})$
is a Markov process on $[0,T]$.

Let $I=\{0,1\}^{n}$. For $\iota=(i_1,\dots,i_n)\in I$ we set
$D_{\iota}=\{x\in\BR^n; (-1)^{i_{k}}x_{k}>0,k=1,\dots,n\}$,
$P=\bigcup_{\iota\in I} D_{\iota}$, $P_{T}=[0,T)\times P$. By
It\^o's formula,
\begin{equation}
\label{eq2.02} X^{s,x,i}_t=x^i\exp\big((r-d_i-a_{ii})(t-s)
+\sum^n_{j=1}\sigma_{ij}(W^j_t-W^j_s)\big),\quad t\in[s,T].
\end{equation}
Therefore if $s\in[0,T)$ and $x\in D_{\iota}$ for some $\iota\in
I$ then $P_{s,x}(X_{t}\in D_{\iota},\, t\ge s)=1$. From this  and
the fact that $a$ is positive definite it follows that if $x\in
P_T$ then $\det\sigma(X_t)>0$, $P_{s,x}$-a.s. for every $t\ge s$,
where $\sigma(x)=\{\sigma_{ij}x_i\}_{i,j=1,\dots, n}$. Moreover,
$[s,T]\ni t\mapsto\sigma^{-1}(X_t)$ is a continuous process.
Therefore, if $x\in P_T$ then by L\'evy's theorem the process
$B_{s,\cdot}$ defined as
$B_{s,t}=\int^t_s\sigma^{-1}(X_{\theta})\,dM_{\theta}$, where
$M^i_t=X^i_t-X^i_0-\int^t_0(r-d_i)X^i_{\theta}\,d\theta$, $t\in
[s,T]$, is under $P_{s,x}$ a standard $n$-dimensional
$\{\FF^s_t\}$-Wiener process on $[s,T]$ and
\begin{equation}
\label{eq2.05}
X^i_t-x^i=\int_{s}^{t}(r-d_{i})X^{i}_{\theta}\,d\theta
+\sum_{j=1}^n\int_{s}^{t}\sigma_{ij}X^{i}_{\theta}\,dB^j_{s,\theta},
\quad t\in[s,T],\quad P_{s,x}\mbox{-a.s.},
\end{equation}
i.e.
\begin{equation}
\label{eq2.04} X^{i}_t=x^i\exp\big((r-d_i-a_{ii})(t-s)
+\sum^n_{j=1}\sigma_{ij}B_{s,t}\big),\quad t\in[s,T],\quad
P_{s,x}\mbox{-a.s.}
\end{equation}
The above forms of the assets price dynamics will be more
convenient for us than (\ref{eq1.1}) or (\ref{eq2.02}). Note that
from the definition of the process $B_{s,\cdot}$ and
(\ref{eq2.04}) it follows that
$\sigma(X_{\theta};\theta\in[s,t])=\sigma(B_{s,\theta};\theta\in[s,t])$
for $s\in[0,T)$, so for every $s\in[0,T)$ the filtration
$\{\FF^s_t\}$ is the completion of the Brownian filtration.

 In Bensoussan \cite{Be} and Karatzas \cite{K}
(see also Section 2.5 in \cite{KS}) it is shown that under (A1)
the arbitrage-free value $V$ of an American option with payoff
function $\psi$ and expiration time $T$ is given by the solution
of the stopping problem
\begin{equation}
\label{eq1.2} V(s,x)=\sup_{\tau\in\mathcal{T}_{s}}
E_{s,x}\big(e^{-r(\tau-s)}\psi(X_{\tau})\big),
\end{equation}
where the supremum is taken over the set $\mathcal{T}_{s}$ of all
$\{\FF^s_t\}$-stopping times $\tau$ with values in $[s,T]$.

From the results proved in  \cite{EKPPQ} it follows  that under
(A1) for every $(s,x)$ there exists a unique solution
$(Y^{s,x},Z^{s,x},K^{s,x})$, on the space
$(\Omega,\FF^s_T,P_{s,x})$, to the reflected BSDE with terminal
condition $\psi(X_T)$, coefficient $f:\BR\rightarrow\BR$ defined
as $f(y)=-ry$, $y\in\BR$, and barrier $\psi(X)$
(RBSDE${}_{s,x}(\psi,-r y,\psi)$ for short). This means that the
processes $Y^{s,x},Z^{s,x},K^{s,x}$ are
$\{\FF^s_t\}$-progressively measurable, satisfy some integrability
conditions and $P_{s,x}$-a.s.,
\begin{equation}
\label{eq2.1} \left\{
\begin{array}{l}
Y^{s,x}_t=\psi(X_T)-\int^T_tr Y^{s,x}_{\theta}\,d\theta
+K^{s,x}_T-K^{s,x}_t
-\int^T_tZ^{s,x}_{\theta}\,dB_{s,\theta},\quad
t\in[s,T],\medskip\\
Y^{s,x}_t\ge \psi(X_t),\quad t\in[s,T],\medskip \\
K^{s,x}\mbox{ is increasing, continuous, }K^{s,x}_s=0,\,\,
\int^T_s(Y^{s,x}_t-\psi(X_t))\,dK^{s,x}_t=0.
\end{array}
\right.
\end{equation}
In \cite{EKPPQ} it is also proved that for every $(s,x)\in Q_T$,
\begin{equation}
\label{eq2.2} Y^{s,x}_t=u(t,X_t),\quad t\in[s,T],\quad
P_{s,x}\mbox{-a.s.},
\end{equation}
where $u$ is a viscosity solution to the obstacle problem
\begin{equation}
\label{eq2.3} \left\{
\begin{array}{ll}
\min(u(s,x)-\psi(x),-u_s-L_{BS}u(s,x)+r u(s,x))=0, &
(s,x)\in Q_T,\medskip\\
u(T,x)=\psi(x), & x\in\BR^n
\end{array}
\right.
\end{equation}
with
\[
L_{BS}u=\sum^n_{i=1}(r-d_{i})x_{i}u_{x_{i}} +\frac12\sum^n_{i,j=1}
a_{ij}x_ix_ju_{x_{i}x_{j}}.
\]
From \cite{EKPPQ,EQ} we  know that $V$ defined by (\ref{eq1.2}) is
equal to $Y^{s,x}_s$. Hence
\begin{equation}
\label{eq2.4}
V(s,x)=Y^{s,x}_s=u(s,x),\quad(s,x)\in[0,T]\times\BR^n.
\end{equation}
In the next section we analyze $V$ via  (\ref{eq2.4}) but as a
matter of fact instead of viscosity  solutions of (\ref{eq2.3}) we
consider variational solutions which provide more information on
the value function $V$.

\section{Obstacle problem for the Black and Scholes equation}
\label{sec3}

Assume that $\psi:\BR^n\rightarrow\BR_+$ is continuous and
satisfies the polynomial growth condition. Let
$L^2_{\varrho}=L^2(\BR^n;\varrho^2\,dx)$, $H^{1}_{\varrho}=\{u\in
L^2_{\varrho}:\sum^n_{j=1} \sigma_{ij}x_iu_{x_j}\in
L^{2}(\BR^n;\varrho^2\,dx),\,i=1,\dots,n\}$ and
$W_{\varrho}=\{u\in L^2(0,T;H^1_{\varrho}):u_t\in
L^2(0,T;H^{-1}_{\varrho})\}$, where $u_t,u_{x_i}$ denote the
partial derivatives in the distribution sense,
$\varrho(x)=(1+|x|^2)^{-\gamma}$ and $\gamma>0$ is chosen so that
$\int_{\BR^n}\varrho^2(x)\,dx<\infty$ and
$\int_{\BR^n}\psi^{2}(x)\varrho^2(x)\,dx<\infty$. Following
\cite{BCFE,KR:BPAN} we adopt the following definition.

\begin{definition}
(a) A pair $(u,\mu)$ consisting of $u\in W_{\varrho}\cap C(\bar
Q_T)$ and a Radon measure $\mu$ on $Q_T$ is a variational solution
to (\ref{eq2.3}) if
\[
u(T,\cdot)=\psi,\quad u\ge
\psi,\quad \int_{Q_T}(u-\psi)\varrho^2\,d\mu=0
\]
and the equation
\[
u_t+L_{BS}u=ru-\mu
\]
is satisfied in the strong sense, i.e. for every $\eta\in
C_{0}^{\infty}(Q_{T})$,
\[
\langle u_t,\eta\rangle_{\varrho,T}+ \langle
L_{BS}u,\eta\rangle_{\varrho,T} =r\langle
u,\eta\rangle_{2,\varrho,T}-\int_{Q_{T}}\eta\varrho^{2}\,d\mu,
\]
where
\[
\langle L_{BS}u,\eta\rangle_{\varrho,T}
=\sum^n_{i=1}\langle(r-d_i)x_iu_{x_i},\eta\rangle_{2,\varrho,T}
-\frac12\sum^n_{i,j=1}a_{ij}\langle u_{x_i},
(x_ix_j\eta\varrho^2)_{x_j}\rangle_{2,T}.
\]
Here $\langle\cdot,\cdot\rangle_{\varrho,T}$ stands for the
duality pairing between $L^2(0,T;H^1_{\varrho})$ and
$L^2(0,T;H^{-1}_{\varrho})$,
$\langle\cdot,\cdot\rangle_{2,\varrho,T}$ is the usual scalar
product in $L^2(0,T;L^2_{\varrho})$ and
$\langle\cdot,\cdot\rangle_{2,T}
=\langle\cdot,\cdot\rangle_{2,\varrho,T}$ with $\varrho\equiv1$.
\medskip\\
(b) If $\mu$ in the above definition admits a density (with
respect to the Lebesgue measure) of the form
$\Phi_u(t,x)=\Phi(t,x,u(t,x))$ for some measurable $\Phi:\bar
Q_T\times\BR\rightarrow\BR_+$, then we say that $u$ is a
variational solution to the semilinear problem
\begin{equation}
\label{eq3.10} u_t+L_{BS}u=ru-\Phi_u,\quad u(T,\cdot)=\psi,\quad
u\ge\psi.
\end{equation}
\end{definition}

In our main theorems below we show that if $\psi$ satisfies (A1)
and (A2) or (A3)  then  the measure $\mu$ is absolutely continuous
with respect to the Lebesgue measure and its density has the form
$\Phi_u(t,x)=\mathbf{1}_{\{u(t,x)=\psi(x)\}}\Psi^{-}(x)$, where
$\Psi^{-}=\max\{-\Psi,0\}$ and $\Psi$ is determined by $\psi$ and
the parameters $r,d,a$. In the next section we compute $\Psi$ for
some concrete options.

\subsection{Payoffs satisfying (A1), (A2)}

\begin{uw}
One can check that if $u$ is a solution to (\ref{eq3.10}) then $v$
defined as
\[
v(t,x)=u(T-t,(-1)^{i_1}e^{x_1},\dots,(-1)^{i_n}e^{x_n})\equiv
u(T-t,e^x)
\]
for $t\in[0,T]$, $x=(x_1,\dots,x_n)\in D_{\iota}$, $\iota\in I$
($D_{\iota}$ is defined in Section \ref{sec2}) is a variational
solution of the Cauchy problem
\[
v_{t}-Lv=-r v+\bar\Phi,\quad v\ge\bar\psi,\quad
v(0,\cdot)=\bar\psi,
\]
where
\[
Lv=\sum^n_{i=1}(r-d_{i}-\frac12 \sigma_{ii}^{2})v_{x_{i}}+
\frac12\sum^n_{i,j=1}a_{ij}v_{x_{i}x_{j}}
\]
and $\bar\Phi(t,x)=\Phi_u(T-t,e^x)$,
$\bar\psi(t,x)=\psi(T-t,e^x)$. Furthermore, a simple calculation
shows that if $\eta$ is a smooth function on $\BR^n$ with compact
support and $U\subset\BR^n$ is a bounded open set such that
supp$[\eta]\subset U$ then $\tilde{v}=v\eta$ is a solution of the
Cauchy-Dirichlet problem
\[
\tilde{v}_{t}-\tilde{L}\tilde{v}=-r\tilde{v}+f,\quad \quad
\tilde{v}(0,\cdot)=\tilde{\psi},\quad \tilde
v_{|[0,T]\times\partial U}=0,
\]
where $\tilde{\psi}=\bar\psi\eta$, $\tilde{L}$ is some uniformly
elliptic operator with smooth coefficients not depending on $t$
and $f\in L^{2}(0,T;L^2(U))$. By classical regularity results
(see, e.g., Theorem 5 in \S7.1 in \cite{E}), $\tilde v\in
L^2(0,T;H^2(U))\cap L^{\infty}(0,T;H^1_0(U))$ and $\tilde v_t\in
L^2(0,T;L^2(U))$. From this and the construction of $\tilde v$ we
infer that the regularity properties of $\tilde v$ are retained by
$u$. It follows in particular that
\begin{equation}
\label{eq3.11} u_t+L_{BS}u=ru-\Phi_u\quad\mbox{a.e. on }P_T.
\end{equation}
\end{uw}

\begin{tw}
\label{tw2.1} Assume \mbox{\rm(A1), (A2)}.
\begin{enumerate}
\item[\rm(i)]$u$ defined by \mbox{\rm(\ref{eq2.4})} is a variational
solution of the semilinear Cauchy problem
\begin{equation}
\label{eq2.9} u_{t}+L_{BS}u=r u-\Phi^{-}_u,\quad u(T,\cdot)=\psi
\end{equation}
with
\[
\Phi_u(t,x)=\mathbf{1}_{\{u(t,x)=\psi(x)\}}\Psi(x),\quad (t,x)\in
Q_T,
\]
where for $x\in\BR^n$ such that $(t,x)\in\{u=\psi\}$,
\[
\Psi(x)=-r\psi(x)+L_{BS}\psi(x).
\]
\item[\rm(ii)]Set $\sigma(x)=\{\sigma_{ij}x_i\}_{i,j=1,\dots,n}$
and
\begin{equation}
\label{eq3.1}
K_{s,t}=\int_{s}^{t}\Phi^{-}_u(\theta,X_{\theta})\,d\theta,\quad
t\in [s,T].
\end{equation}
Then for every $(s,x)\in P_{T}$ the triple
$(u(\cdot,X),\sigma(X)u_x(\cdot,X),K_{s,\cdot})$ is a unique
solution of RBSDE$_{s,x}(\psi,-r y,\psi)$.
\end{enumerate}
\end{tw}
\begin{dow}
Fix $(s,x)\in P_{T}$. Let $(Y^{s,x},Z^{s,x},K^{s,x})$ be a
solution of RBSDE$_{s,x}(\psi,-r y,\psi)$ and let $u$ be a
viscosity solution of (\ref{eq2.3}). For $t_0\in(s,T)$ let
$U\subset\BR^n$ be an open set of assumption (A2). Then there
exists $\eta\in C^{\infty}(\BR^n)$ such that $\eta\ge0$,
$\eta\equiv1$ on $\{u=\psi\}\cap Q_{t_0}$ and $\eta\equiv0$ on
$U^{c}$ (we make the convention that $\eta(t,x)=\eta(x)$). Of
course $(Y^{s,x},Z^{s,x}, K^{s,x})$ is a solution of
RBSDE$_{s,x}(Y^{s,x}_{t_{0}},-r y,\psi)$ on $[s,t_{0}]$. It is
also a solution of RBSDE$_{s,x}(Y^{s,x}_{t_{0}},-ry,
\tilde{\psi})$ on $[s,t_{0}]$ with
$\tilde{\psi}(x)=\eta(x)\psi(x)$, because $\tilde\psi\le\psi$ and
by (\ref{eq2.1}) and (\ref{eq2.2}),
\begin{align*}
\int^{t_0}_s(Y^{s,x}_t-\tilde\psi(X_t))\,dK^{s,x}_t
&=\int^{t_0}_s(u(t,X_t)-\tilde\psi(X_t))
\mathbf{1}_{\{u(t,X_t)=\psi(X_t)\}}\,dK^{s,x}_t\\
&=\int^{t_0}_s(u(t,X_t)-\psi(X_t))\,dK^{s,x}_t=0.
\end{align*}
Since $\tilde{\psi}$ is smooth, applying It\^o's formula yields
\[
\tilde{\psi}(X_{t}) =\tilde\psi(X_s)+
\sum^n_{i=1}\int_{s}^{t}\tilde{\psi}_{x_i}(X_{\theta})\,dX^i_{\theta}
+\frac12\sum^n_{i,j=1}\int_{s}^{t}a_{ij}X^i_{\theta}X^j_{\theta}
\tilde{\psi}_{x_{i}x_{j}}(X_{\theta})\,d\theta.
\]
From the above, (\ref{eq2.2}) and \cite[Remark 4.3]{EKPPQ} it
follows that there exists a predictable process $\alpha^{s,x}$
such that $0\le \alpha^{s,x}\le 1$ and
\begin{align*}
dK^{s,x}_{t}=\alpha^{s,x}_{t}\mathbf{1}_{\{u={\psi}\}}(X_{t})
\Big(-r\tilde{\psi}(X_t)
&+\sum^n_{i=1}(r-d_{i})X^{i}_{t}\tilde{\psi}_{x_{i}}(X_{t})\\
&\qquad+\frac12\sum^n_{i,j=1} a_{ij}
X^i_tX^j_t\tilde{\psi}_{x_{i}x_{j}}(X_{t})\Big)^{-}\,dt
\end{align*}
on $[s,t_{0}]$. Thus
\begin{equation}
\label{eq3.20}
dK^{s,x}_{t}=\alpha^{s,x}_{t}\mathbf{1}_{\{u(t,X_t)=\psi(X_t)\}}
\Psi^{-}(X_t)\,dt
\end{equation}
on $[s,t_{0}]$ for every $t_0\in[s,T)$. Consequently, the above
equation is satisfied on $[s,T]$. Since the coefficients of the
stochastic differential equation (\ref{eq2.05}) satisfy the
assumptions of the ``equivalence of norm" result proved in
\cite{BM} (see \cite[Proposition 5.1]{BM}), it follows from
\cite[Theorem 3]{BCFE} that there exists a function $\alpha$ on
$Q_{T}$ such that $0\le\alpha\le 1$ a.e. and for a.e. $(s,x)\in
Q_{T}$,
\begin{equation}
\label{eq3.2} \alpha^{s,x}_{t}=\alpha(t,X_{t}),\,\quad dt\otimes
P_{s,x}\mbox{-a.s.}
\end{equation}
Moreover, $u\in C(\bar Q_T)$ by \cite[Lemma 8.4]{EKPPQ} and from
\cite[Theorem 3]{BCFE} it follows that $u\in W_{\varrho}$ and $u$
is a variational solution of the Cauchy problem
\[
u_{t}+L_{BS}u=r u-\alpha\mathbf{1}_{\{u=\psi\}}\Psi^{-},\quad
u(T,\cdot)=\psi.
\]
By the above and (\ref{eq3.11}),
\[
u_{t}+L_{BS}u=r u -\alpha\mathbf{1}_{\{u=\psi\}}\Psi^{-}\quad
\mbox{a.e. on } Q_T,
\]
so by Lemma A.4 in Chapter II in \cite{KiS},
\[
\psi_{t}+L_{BS}\psi=r\psi-\alpha\Psi^{-}\quad \mbox{a.e. on }
\{u=\psi\}.
\]
On the other hand, by the definition of $\Psi$,
\[
\psi_{t}+L_{BS}\psi=L_{BS}\psi=r\psi+\Psi\quad \mbox{on }
\{u=\psi\}.
\]
Thus $\Psi=-\alpha\Psi^{-}$ a.e. on $\{u=\psi\}$, which implies
that $\alpha\Psi=\Psi$ a.e. on $\{u=\psi\}$,  and hence that
\begin{equation}
\label{eq2.10}
\mathbf{1}_{\{u=\psi\}}\alpha\Psi^{-}
=\mathbf{1}_{\{u=\psi\}}\Psi^{-}\quad\mbox{a.e.}
\end{equation}
Accordingly (\ref{eq2.9}) is satisfied. From (\ref{eq2.02}) it is
clear that if $s\in[0,T)$ and $x\in D_{\iota}$ for some $\iota\in
I$ then $P_{s,x}(X_{t}\in D_{\iota},\, t\ge s)=1$ and for every
$t\in(s,T]$ the random variable $X_t$ has strictly positive
density on $D_{\iota}$ under $P_{s,x}$. From this and
(\ref{eq2.10}) it follows that
\begin{equation}
\label{eq3.3}
\mathbf{1}_{\{u=\psi\}}(t,X_{t})\alpha(t,X_{t})\Psi^{-}(X_t)
=\mathbf{1}_{\{u=\psi\}}(t,X_{t})\Psi^{-}(X_t),\quad dt\otimes
P_{s,x}\mbox{-a.s.}
\end{equation}
for every $(s,x)\in P_{T}$. In \cite{Kl:PA1} it is proved that the
function $\mathbf{1}_{\{u=\psi\}}\alpha$ is a weak limit in
$L^{2}(Q_{T})$ of some sequence $\{\alpha_{n}\}$ of nonnegative
functions bounded by $1$ and such that
$\alpha_{n}(t,X_{t})\rightarrow \alpha^{s,x}_{t}$ weakly in
$L^{2}([0,T]\times\Omega;dt\otimes P_{s,x})$ for every $(s,x)\in
Q_{T}$. Therefore using once again the fact that for every
$(s,x)\in P_{T}$ the process $X$ has a strictly positive
transition density under $P_{s,x}$ we conclude that (\ref{eq3.2})
holds for every $(s,x)\in P_{T}$, which when combined with
(\ref{eq3.3}) implies (\ref{eq3.1}). What is left is to prove that
for every $(s,x)\in P_{T}$,
\begin{equation}
\label{eq3.4} Z^{s,x}_{t}=\sigma(X_t)u_x (t,X_{t}),\quad dt\otimes
P_{s,x}\mbox{-a.s.}
\end{equation}
From the results proved in \cite[Section 6]{EKPPQ} it follows that
for every $(s,x)\in Q_{T}$,
\begin{equation}
\label{eq3.5} E_{s,x}\sup_{s\le t\le T}
|Y^{s,x,n}_t-Y^{s,x}_{t}|^{2}
+E_{s,x}\int_{s}^{T}|Z^{s,x,n}_{t}-Z^{s,x}_{t}|^{2}\,dt\rightarrow0,
\end{equation}
where $(Y^{s,x,n},Z^{s,x,n})$ is a solution of the  BSDE
\begin{align}
\label{eq2.18} Y^{s,x,n}_{t}&=\psi(X_{T})-\int_{t}^{T}r
Y^{s,x,n}_{\theta}\,d\theta \nonumber\\
&\quad+\int_{t}^{T}n(Y^{s,x,n}_{\theta}-\psi(X_{\theta}))^{-}\,d\theta
-\int_{t}^{T}Z^{s,x,n}_{\theta}\,dB_{s,\theta}.
\end{align}
It is known (see \cite{Pardoux}) that
\begin{equation}
\label{eq3.8} Y^{s,x,n}_{t}=u_{n}(t,X_{t}),\quad t\in [s,T],\quad
P_{s,x}\mbox{-a.s.},
\end{equation}
where $u_{n}$ is a viscosity solution of the Cauchy problem
\[
(u_n)_t+L_{BS}u_n=-r u_n+n(u_n-\psi)^{-},\quad u_n(T,\cdot)=\psi.
\]
We know that $P_{s,x}(X_{t}\in D_{\iota},\, t\ge s)=1$ if $x\in
D_{\iota}$. Moreover, by classical regularity results (see, e.g.,
\cite[Theorem 1.5.9]{F} and Remark preceding Theorem \ref{tw2.1}),
$u_{n}\in C^{1,2}(P_{T})$. Therefore applying It\^o's formula
shows that (\ref{eq2.18}) holds true with $Z^{s,x,n}_{\theta}$
replaced by $\sigma(X_{\theta})(u_n)_x(\theta,X_{\theta})$. Since
(\ref{eq2.18}) has a unique solution (see \cite[Corollary
3.7]{EKPPQ}), it follows that
\begin{equation}
\label{eq3.29} Z^{s,x,n}_{t}=\sigma(X_t)(u_{n})_x(t,X_{t}), \quad
dt\otimes P_{s,x}\mbox{-a.s.}
\end{equation}
for every $(s,x)\in P_{T}$. By (\ref{eq3.5}) and (\ref{eq3.8}),
$u_n\rightarrow u$ pointwise in $Q_T$. Moreover, from
(\ref{eq3.8}), (\ref{eq3.29}) and standard estimates for solutions
of BSDEs (see, e.g., \cite[Section 6]{EKPPQ}) it follows that
there is $C>0$ such that for any $(s,x)\in P_T$,
\begin{equation}
\label{eq3.30} E_{s,x}\sup_{s\le t\le T}
|u_n(t,X_t)|^2+E_{s,x}\int^T_s|\sigma(X_t)(u_n)_x(t,X_t)|^2\,dt
\le CE_{s,x}|\psi(X_T)|^2,
\end{equation}
while from (\ref{eq3.5}), (\ref{eq3.29}) it follows that
\begin{equation}
\label{eq3.31} E_{s,x}\sup_{s\le t\le T}
\int^T_s|\sigma(X_t)((u_n)_x-(u_m)_x)(t,X_t)|^2\,dt \rightarrow0
\end{equation}
as $n,m\rightarrow\infty$. From (\ref{eq3.30}) one can deduce that
$u_n\in L^2(0,T;H_{\varrho})$ and then, by using (\ref{eq3.31}),
that $u_n\rightarrow u$ in $L^2(0,T;H_{\varrho})$ (see the
arguments following (2.12) in the proof of \cite[Theorem
2.3]{KR:BPAN}). From the last convergence and (\ref{eq3.5}),
(\ref{eq3.29}) it may be concluded that
\[
E_{s,x}\int^T_s|\sigma(X_t)(u_n)_x(t,X_t)-Z^{s,x}_t|^2\,dt=0
\]
for $(s,x)\in P_T$, which implies (\ref{eq3.4}).
\end{dow}

\subsection{Convex payoffs}

Assume that $\psi:\BR^n\rightarrow\BR$ is convex. Let $m$ denote
the Lebesgue measure on $\BR^n$, $\nabla_i\psi$ denote the usual
partial derivative with respect to $x_i$, $i=1\dots,n$, and let
$E$ be set of all $x\in\BR^n$ for which the gradient
\[
\nabla\psi(x)=(\nabla_1\psi(x),\dots,\nabla_n\psi(x))
\]
exists. Since $\psi$ is locally Lipschitz function, $m(E^c)=0$ and
$\nabla\psi=(\psi_{x_1},\dots,\psi_{x_n})$ a.e. (recall that
$\psi_{x_i}$ stands for the partial derivative in the distribution
sense). Moreover, for a.e. $x\in E$ there exists an
$n$-dimensional symmetric matrix $\{H(x)=\{H_{ij}(x)\}$ such that
\begin{equation}
\label{eq3.26} \lim_{E\ni y\rightarrow x}
\frac{\nabla\psi(y)-\nabla\psi(x) -H(x)(y-x)}{|y-x|}=0,
\end{equation}
i.e. $H_{ij}(x)$ are defined as limits through the set where
$\nabla_i\psi$ exists (see, e.g., \cite[Section 7.9]{AA}). By
Alexandrov's theorem (see, e.g., \cite[Theorem 7.10]{AA}), if
$x\in E$ is a point where (\ref{eq3.26}) holds then $\psi$ has
second order differential at $x$ and $H(x)$ is the hessian matrix
of $\psi$ at $x$, i.e. $H(x)=\{\nabla^2_{ij}\psi(x)\}$.

The second order derivative of $\psi$ in the distribution sense
$D^2\psi=\{\psi_{x_ix_j}\}$ is a matrix of real-valued Radon
measures $\{\mu_{ij}\}$ on $\BR^n$ such that $\mu_{ij}=\mu_{ji}$
and for each Borel set $B$, $\{\mu_{ij}(B)\}$ is a nonnegative
definite matrix (see, e.g., \cite[Section 6.3]{EG}). Let
$\mu_{ij}=\mu^a_{ij}+\mu^s_{ij}$ be the Lebesgue decomposition of
$\mu_{ij}$ into the  absolutely continuous and singular parts with
respect to $m$. By Theorem 1 in Section 6.4 in \cite{EG},
\begin{equation}
\label{eq3.17} \mu^a_{ij}(dx)=\nabla^2_{ij}\psi(x)\,dx.
\end{equation}

For $R>0$ set $D_R=P\cap\{x\in\BR^n:|x|<R\}$ and
$\tau_R=\inf\{t\ge s:X_t\notin D_R\}$. Let $\tilde L_{BS}$ denote
the operator formally adjoint to $L_{BS}$. By \cite[Theorem
4.2.5]{Pi} for a sufficiently large $\alpha>0$ there exist the
Green's functions $G^{\alpha}_R$, $\tilde G^{\alpha}_R$ for
$\alpha-L_{BS}$ and $\alpha-\tilde L_{BS}$ on $D_R$. Let $A$ be a
continuous additive functional of $\BX$ and let $\nu$ denote the
Revuz measure of $A$ (see, e.g., \cite{Re}). By the theorem proved
in Section V.5 of \cite{Re}, for every  nonnegative $f\in
C_0(\BR^d)$,
\[
E_{s,x}\int^{\tau_R}_se^{-\alpha t}f(X_t)\,dA^{\nu}_t
=\int_{\BR^n}G^{\alpha}_R(x,y)f(y)\,\nu(dy).
\]
Since $G^{\alpha}_R(x,y)=\tilde G^{\alpha}_R(y,x)$ by
\cite[Corollary 4.2.6]{Pi}, it follows that
\begin{equation}
\label{eq3.14} E_{s,g\cdot m}\int^{\tau_R}_se^{-\alpha t}
f(X_t)\,dA^{\nu}_t=\int_{\BR^n}\tilde
G^{\alpha}_Rg(y)f(y)\,\nu(dy)
\end{equation}
for any nonnegative $g\in C_0(D_R)$, where $E_{s,g\cdot m}$
denotes the expectation with respect to the measure $P_{s,g\cdot
m}(\cdot)=\int P_{s,x}(\cdot)g(x)\,dx$ and
\[
\tilde G^{\alpha}_Rg(y)=\int_{G_R}\tilde
G^{\alpha}_R(y,x)g(x)\,dx.
\]
Note that if $g$ is not identically equal to zero then $\tilde
G^{\alpha}_Rg$ is strictly positive (see \cite[Theorem
4.2.5]{Pi}).

Set
\[
\LL_{BS}=\sum^n_{i=1}(r-d_{i})x_{i}\nabla_i +\frac12\sum^n_{i,j=1}
a_{ij}x_ix_j\nabla^2_{ij}\,.
\]
\begin{tw}
\label{tw3.2} Assume \mbox{\rm(A3)}. Then assertions \mbox{\rm(i),
(ii)} of Theorem \ref{tw2.1} hold true with $L_{BS}$ replaced by
$\LL_{BS}$.
\end{tw}
\begin{dow}
We use the notation of Theorem \ref{tw2.1}. Fix $s\in[0,T)$. Since
$\psi$ is a continuous convex function, from It\^o's formula
proved in \cite{Bo} it follows that  there exists a continuous
increasing process $A$ such that for $x\in\BR^n$,
\begin{equation}
\label{eq3.9} \psi(X_{t})=\psi(X_{s})+A_{t}
+\int_{s}^{t}\nabla\psi (X_{\theta})\, dX_{\theta},\quad t\in
[s,T],\quad P_{s,x}\mbox{-a.s.}
\end{equation}
From (\ref{eq3.9}) it follows that $A$ is a positive continuous
additive functional (PCAF for short) of $\BX$. Let $\nu$ denote
the Revuz measure of $A$. We are going to show that
$\mathbf{1}_P\cdot\nu=\mathbf{1}_P\cdot\mu $ where $\mu$ is the
measure on $\BR^n$ defined as
\[
\mu(dx)=\sum^n_{i,j=1}a_{ij}x_ix_j\,\mu_{ij}(dx).
\]
To this end, let us set
\[
\mu^{\varepsilon}_{ij}=\frac{\partial^2\psi}{\partial x_i\partial
x_j}*\rho_{\varepsilon}\,,\quad \mu^{\varepsilon}(dx)
=\sum^n_{i,j=1}a_{ij}x_ix_j\,\mu^{\varepsilon}_{ij}(dx),
\]
where $\{\rho_{\varepsilon}\}_{\varepsilon>0}$ is some family of
mollifiers. Fix a nonnegative $g\in C_0(D_R)$  such that $g(x)>0$
for some $x\in D_R$ and denote by $A^{\varepsilon}$ the PCAF of
$\BX$ in Revuz correspondence with $\mu_{\varepsilon}$. Then for a
sufficiently large $\alpha>0$,
\begin{equation}
\label{eq3.15} E_{s,g\cdot m}\int^{\tau_R}_se^{-\alpha t}
f(X_t)\,dA^{\varepsilon}_t =\int_{\BR^n}\tilde
G^{\alpha}_Rg(y)f(y)\,\mu^{\varepsilon}(dy)
\end{equation}
for all nonnegative $f\in C_0(\BR^d)$. By \cite[Theorem 2]{CP},
$E_{s,x}\sup_{t\ge s}
|A^{\varepsilon}_{t\wedge\tau_R}-A_{t\wedge\tau_R}|\rightarrow0$
as $\varepsilon\downarrow0$  for every $x\in\BR^d$. Hence
$\int^{\tau_R}_se^{-\alpha t}f(X_t)\,dA^{\varepsilon}_t
\rightarrow \int^{\tau_R}_se^{-\alpha t}f(X_t)\,dA_t$ weakly under
$P_{s,x}$ for $x\in\BR^d$. Since
\[
A^{\varepsilon}_{t\wedge\tau_R}=\psi_{\varepsilon}(X_{t\wedge\tau_R})
-\psi_{\varepsilon}(X_s)-\int^{t\wedge\tau_R}_s
\nabla\psi_{\varepsilon}(X_{\theta})\,dX_{\theta}
\]
and $\sup_{\varepsilon>0}\sup_{|x|\le R}|
\nabla\psi_{\varepsilon}(x)|\le C(R)<\infty$ by Lemma in
\cite{CP}, it follows that for every compact subset
$K\subset\BR^n$, $\sup_{x\in K}
\sup_{\varepsilon>0}E_{s,x}|A^{\varepsilon}_{t\wedge\tau_R}|^2<\infty$.
Therefore
\begin{equation}
\label{eq3.16} E_{s,g\cdot m}\int^{\tau_R}_se^{-\alpha t}
f(X_t)\,dA^{\varepsilon}_t\rightarrow E_{s,g\cdot m}
\int^{\tau_R}_se^{-\alpha t}f(X_t)\,dA_t
\end{equation}
as $\varepsilon\downarrow0$. On the other hand, since
$\mu^{\varepsilon}_{ij}\rightarrow\mu_{ij}$ weakly${}^{*}$ for
$i,j=1,\dots,n$ and, by \cite[Theorem 4.2.5]{Pi}, $f\tilde
G^{\alpha}_Rg\in C_0(G_R)$, we have
\[
\sum^n_{i,j=1}\int_{\BR^n} \tilde
G^{\alpha}_Rg(y)f(y)a_{ij}y_iy_j\, \mu^{\varepsilon}_{ij}(dy)
\rightarrow \sum^n_{i,j=1}\int_{\BR^n} \tilde
G^{\alpha}_Rg(y)f(y)\, a_{ij}y_iy_j\mu_{ij}(dy).
\]
Combining this with (\ref{eq3.14}), (\ref{eq3.15}), (\ref{eq3.16})
we see that for every $f\in C_0(\BR^n)$,
\[
\int_{\BR^n}\tilde G^{\alpha}_Rg(y)f(y)\mu(dy)=\int_{\BR^n}\tilde
G^{\alpha}_Rg(y)f(y)\,\nu(dy).
\]
Since $\tilde G^{\alpha}_Rg$ is strictly positive on $D_R$, we
conclude from the above that $\mu=\nu$ on $D_R$ for each $R>0$.
Consequently, $\mu=\nu$ on $P$. For $x\in P$,
$P_{s,x}(X_t\in\BR^n\setminus P)=0$ for $t\ge s$. Hence
\begin{equation}
\label{eq3.18}
A^{\nu}_t=\int^t_s\mathbf{1}_P(X_s)\,dA^{\nu}_s=A^{\mathbf{1}_P\cdot\nu}_t
=A^{\mathbf{1}_P\cdot\mu}_t,\quad t\ge s,\quad P_{s,x}\mbox{-a.s.}
\end{equation}
for $x\in P$. Let $\mu^a$ denote the absolutely continuous part in
the Lebesgue decomposition of $\mathbf{1}_P\cdot\mu$. By
(\ref{eq3.17}),
$\mu^a(dx)=\sum^n_{i,j=1}\mathbf{1}_P(x)a_{ij}x_ix_j\nabla^2_{ij}\psi(x)\,dx$.
Hence
\begin{equation}
\label{eq3.19}
A^{\mu^a}_t=\sum^n_{i,j=1}\int^t_sa_{ij}X^i_{\theta}X^j_{\theta}
\nabla^2_{ij}\psi(X_{\theta})\,d\theta, \quad t\ge s, \quad
P_{s,x}\mbox{-a.s.}
\end{equation}
for $x\in P$. From (\ref{eq3.9}), (\ref{eq3.18}), (\ref{eq3.19})
and \cite[Remark 4.3]{EKPPQ} it follows that
\begin{align*}
dK^{s,x}_{t}=\alpha^{s,x}_{t}\mathbf{1}_{\{u=\psi\}}
(X_{t})\Big(-r\psi(X_{t})
&+\sum_{i=1}^{n}(r-d_{i})X^{i}_{t}\nabla_i\psi(X_{t})\\
&\qquad+\frac12\sum^n_{i,j=1} a_{ij}
X^i_tX^j_t\nabla^2_{ij}\psi(X_{t})\Big)^{-}\,dt.
\end{align*}
Let $u$ be a viscosity solution of (\ref{eq2.3}). From the above
and the results proved in \cite{BCFE} (see the reasoning following
(\ref{eq3.20})) we conclude that $u\in W_{\varrho}\cap C(\bar
Q_{T})$ and there is a function $\alpha$ on $Q_T$ such that
$0\le\alpha\le 1$ a.e., (\ref{eq3.2}) is satisfied and $u$ is a
variational solution of the Cauchy problem
\begin{equation}
\label{eq3.24} u_{t}+L_{BS}u=r u
-\alpha\mathbf{1}_{\{u=\psi\}}\Psi^{-},\quad u(T,\cdot)=\psi
\end{equation}
with
\begin{equation}
\label{eq3.22} \Psi=-r\psi+\LL_{BS}\psi\quad \mbox{on }
\{u=\psi\}.
\end{equation}
By Remark preceding Theorem \ref{tw2.1}, $u(t,\cdot)\in
H^2_{loc}(\BR^n)$. Therefore by Remark (ii) following Theorem 4 in
Section 6.1 in \cite{EG} the distributional derivatives $u_{x_i}$,
$u_{x_ix_j}$ are a.e. equal to the approximate derivatives
$\nabla^{ap}_iu$, $(\nabla^{ap})^2_{ij}u$. Let $\LL^{ap}_{BS}$
denote the operator defined as $\LL_{BS}$ but with $\nabla_i$,
$\nabla_{ij}$ replaced by $\nabla^{ap}_i$, $(\nabla^{ap})^2_{ij}$.
Then $u$ is a variational solution of (\ref{eq3.24}) with $L_{BS}$
replaced by $\LL^{ap}_{BS}$ and (\ref{eq3.11}) holds with $L_{BS}$
replaced by $\LL^{ap}_{BS}$.
Hence
\[
u_{t}+\LL^{ap}_{BS}u=ru
-\alpha\mathbf{1}_{\{u=\psi\}}\Psi^{-}\quad\mbox{ a.e. on }Q_T.
\]
On the other hand, since $\psi$ is convex, $\psi\in
BV_{loc}(\BR^n)$ as a locally Lipschitz continuous function and,
by Theorem 3 in Section 6.3 in \cite{EG}, $\psi_{x_i}\in
BV_{loc}(\BR^n)$, $i=1,\dots,n$. Therefore $\psi$ is twice
approximately differentiable a.e. by Theorem 4 in Section 6.1 in
\cite{EG}. It follows now from Theorem 3 in Section 6.1 in
\cite{EG} that $\LL^{ap}u=\LL^{ap}\psi$ a.e. on $\{u=\psi\}$.
Consequently,
\begin{equation}
\label{eq3.23}
\LL^{ap}_{BS}\psi=r\psi-\alpha\Psi^{-}\quad\mbox{a.e. on }
\{u=\psi\}.
\end{equation}
Moreover, since $\psi$ is convex, $\LL_{BS}\psi=\LL^{ap}_{BS}\psi$
a.e. on $\BR^n$ by Remark (i) following Theorem 4 in Section 6.1
in \cite{EG}. 
Therefore combining (\ref{eq3.22}) with (\ref{eq3.23}) we see that
$\Psi=-\alpha\Psi^{-}$ a.e. on $\{u=\psi\}$ from which as in the
proof of Theorem \ref{tw2.1} we get (\ref{eq3.3}). To complete the
proof it suffices now to repeat step by step the arguments
following (\ref{eq3.3}) in the proof of Theorem \ref{tw2.1}.
\end{dow}

\section{The early exercise premium representation}

Let $\xi$ denote the payoff process for an American option with
payoff function $\psi$, i.e.
\[
\xi_t=e^{-r(t-s)}\psi(X_t),\quad t\in[s,T],
\]
and let $\eta$ denote the Snell envelope for $\xi$, i.e. the
smallest supermartingale which dominates $\xi$. It is known (see,
e.g., Section 2.5 in \cite{KS}) that
\[
\eta_t=e^{-r(t-s)}V(t,X_t),\quad t\in[s,T].
\]
Assume (A1), (A2) or (A3). Applying It\^o's formula and using
Theorem \ref{tw2.1} or \ref{tw3.2} we get
\begin{align*}
\eta_t=e^{-r(t-s)}Y^{s,x}_t&=e^{-r(T-s)}\psi(X_T) +\int^T_t
e^{-r(\theta-s)}\Phi^{-}(X_{\theta},Y^{s,x}_{\theta})\,d\theta\\
&\quad-\int^T_te^{-r(\theta-s)}Z^{s,x}_{\theta}\,dW_{\theta},\quad
t\in[s,T],\quad P_{s,x}\mbox{-a.s.},
\end{align*}
which leads to the following corollary.
\begin{corollary}
For every $(s,x)\in Q_T$ the  Snell envelope admits the
representation
\begin{equation}
\label{eq4.1} \eta_t=E_{s,x}\Big(e^{-r(T-s)}\psi(X_T)
+\int^T_te^{-r(\theta-s)}
\Phi^{-}(X_{\theta},Y^{s,x}_{\theta})\,d\theta\,|\FF_t\Big),\quad
t\in[s,T].
\end{equation}
\end{corollary}

Taking $t=s$ in (\ref{eq4.1}) and using (\ref{eq2.2}) we get the
early exercise premium representation for the value function.

\begin{corollary}
For every $(s,x)\in Q_T$ the value function $V$ admits the
representation
\[
V(s,x)=V^E(s,x)+E_{s,x}\int^T_s e^{-r(t-s)}
{\mathbf{1}}_{\{V(t,X_t)=\psi(X_t)\}}\Psi^{-}(X_t)\,dt,
\]
where
\[
V^E(s,x)=E_{s,x}\big(e^{-r(T-s)}\psi(X_T)\big)
\]
is the value of the European option with payoff function $\psi$
and expiration time $T$.
\end{corollary}

In closing this section we show by examples that for many options
$\Psi^{-}$ can be explicitly computed. Using results of \S4 and
\S5 in \cite{Ro} one can check that the payoff functions $\psi$ in
examples 1--4 below satisfy (A3). It is also easy to see that the
payoff function $\psi$ in example 5 satisfies (A2$'$). Note that
the payoff function in example 1 also satisfies (A2$'$) and, by
\cite{BD2,LS}, the payoff functions in examples 2--4 satisfy (A2).
We would like to stress that the last assertion is by no means
evident. On the other hand, the convexity of $\psi$ in examples
2--4 is readily checked.

In all the examples we have computed the corresponding functions
$\Psi^{-}$ on the region $\{u=\psi\}$. When computing $\Psi$ we
keep in mind that $\{u=\psi\}\subset[0,T]\times\{\psi>0\}$.
\begin{enumerate}
\item \underline{Index options and spread options}
\[
\psi(x)=\big(\sum_{i=1}^{n} w_{i}x_i-K\big)^{+},\quad
\Psi^-(x)=\big(\sum_{i=1}^{n} w_{i}d_{i}x_i-r K\big)^{+}\quad
\mbox{(call)}
\]
\[
\psi(x)=\big(K-\sum_{i=1}^{n} w_{i}x_i\big)^{+},\quad
\Psi^-(x)=\big(r K-\sum_{i=1}^{n} w_{i}d_{i}x_i\big)^{+}\quad
\mbox{(put)}
\]
(Here $w_i\in\BR$ for $i=1,\dots,n$).

\item \underline{Max options}
\[
\psi(x)=(\max\{x_1,\dots,x_n\}-K)^{+} \quad \mbox{(call on max)}
\]
\[
\Psi^{-}(x)=\big(\sum_{i=1}^{n} d_{i}\mathbf{1}_{B_{i}}(x)x_i-r
K\big)^{+},
\]
where $B_{i}=\{x\in\BR^n; x_{i}>x_{j},\, j\neq i\}$.

\item \underline{Min options}
\[
\psi(x)=(K-\min\{x_1,\dots,x_n\})^{+}\quad \mbox{(put on min)}
\]
\[
\Psi^{-}(x)=\big(r K
-\sum_{i=1}^{n}d_{i}\mathbf{1}_{C_{i}}(x)x_i)^{+},
\]
where $C_{i}=\{x\in\BR^n; x_{i}<x_{j},\, j\neq i\}$.

\item \underline{Multiple strike options}

\[
\psi(x)=(\max\{x_1-K_{1},\dots, x_n-K_{n}\})^{+},
\]
\[
\Psi^{-}(x)=\big(\sum_{i=1}^{n} \mathbf{1}_{B_{i}}(x-K)(d_{i}x_i-r
K_{i})\big)^{+},
\]
where $K=(K_{1},\dots,K_{n})$.

\item \underline{Power-product options}
\[
\psi(x)=(|x_1\cdot\ldots\cdot x_n|^{\gamma}-K)^+\quad \mbox{for
some }\gamma>0.
\]
If $x\in D_{\iota}$ with $\iota=(i_1,\dots,i_n)\in\{0,1\}^n$ then
\[
\Psi^{-}(x)=\big((r-\gamma\sum_{i=1}^n(r-d_i-a_{ii})
-\gamma^2\sum^n_{i,j=1}a_{ij})f(x)-r K\big)^{+},
\]
where $f(x)=((-1)^{|\iota|}\,x_1\cdot\ldots\cdot x_n)^{\gamma}$
and $|\iota|=i_1+\ldots+ i_n$.
\end{enumerate}

\end{document}